\documentclass[11pt,twoside]{amsart}
\usepackage[all]{xy}
\usepackage {amssymb,latexsym,amsthm,amsmath,mathtools}
\usepackage{enumerate,color}

\usepackage[colorlinks, linkcolor=blue, citecolor=blue]{hyperref}
\usepackage{mathtools}
\usepackage{cleveref}

\long\def\symbolfootnote[#1]#2{\begingroup%
	\def\thefootnote{\fnsymbol{footnote}}\footnote[#1]{#2}\endgroup}

\newcommand{\GL}{\textup{GL}}
\newcommand{\SL}{\textup{SL}}

\makeatletter
\def\imod#1{\allowbreak\mkern10mu({\operator@font mod}\,\,#1)}
\makeatother

\newcommand{\qq}{\mathbb Q}

\newcommand{\nat}{\mathbb{N}}

\newcommand{\mg}[1]{{#1}^{\times}}

\newcommand{\br}[1]{{\rm Br}(#1)}

\newcommand{\trdeg}{\mathrm{tr.deg}}
\newcommand{\trace}{\mathrm{tr}}

\newcommand{\dgal}{\mathsf{DGal}}
\newcommand{\ddeg}{\mathsf{Ddeg}}
\newcommand{\Der}{\mathsf{Der}}
\newtheorem{theorem}{Theorem}[section]

\newtheorem{proposition}[theorem]{Proposition}

\newtheorem{question}[theorem]{Question}
\newtheorem*{theorem*}{Theorem}
\theoremstyle{definition}

\newtheorem{example}[theorem]{Example}
\numberwithin{equation}{section}

\newcommand{\ignore}[1]{}

\newcommand{\mynote}[1]{}

\begin{document}
	\setcounter{section}{0}
	\title{Differential Central Simple Algebras}
	\author{Parul Gupta}
	\email{parulgupta1211@gmail.com}
	\author{Yashpreet Kaur}
	\email{yashpreetkm@gmail.com }
	\author{Anupam Singh}
	\email{anupamk18@gmail.com}
	\address{IISER Pune, Dr. Homi Bhabha Road, Pashan, Pune 411 008, India}
	\thanks{The second named author is supported by NBHM (DAE, Govt. of India): 0204/16(33)/2020/R$\&$D-II/26. The third named author is funded by SERB through CRG/2019/000271 for this research.}
	\date{}
	
	\begin{abstract}
		Differential central simple algebras are the main object of study in this survey article. We recall some crucial notions such as differential subfields, differential splitting fields, tensor products etc. Our main focus is on differential splitting fields which connects these objects to the classical differential Galois theory. 
		We mention several known results and raise some questions along the line.     
	\end{abstract}
	
	\maketitle
	
	\section{Introduction}

	Let $R$ be a ring with identity (not necessarily commutative). A ring $R$ equipped with an additive map $d_R:R\rightarrow R$ satisfying the Leibnitz rule  $d_R(\alpha\beta)=\alpha d_R(\beta)+d_R(\alpha)\beta$ for all $\alpha,\beta \in R$, is called a {\bf differential ring}, and the map $d_R$ is called a {\bf derivation}. We denote a differential ring by an ordered pair $(R,d_R)$. Any ring can be made into a differential ring in a trivial way by taking the zero derivation on it. A differential ring which is a field (or algebra) will be called a {\bf differential field (or algebra)}.  Here are some examples.

	\begin{example}
		The zero derivation is the only derivation on $\qq$. 
	\end{example}

	\begin{example}
		Let  $k =\qq (t)$ be the rational function field over $\qq$.
		For $a\in k$, setting $d_k(t) = a $ we get a derivation on $k$.
		For $a=1$, we get the usual derivation $\frac{d}{dt}$.
	\end{example}

	\begin{example}[Inner derivations]
		Let $R$ be a ring. For any element $a \in R$, the map $\partial_a: R\rightarrow R $ given by $\partial_a (x)= xa -ax$ for $x\in R$ is a derivation on $R$. We say that $\partial_a$ is the \textbf{inner derivation} with respect to $a$. 
		Note that $\partial_a$ restricts to the zero derivation on the center of $R$. In particular, if $R$ is commutative then $\partial_a$ is essentially the zero derivation on $R$ for any $a\in R$. 
		Inner derivations were introduced and studied in the context of Lie algebras in \cite{Jac}. 
	\end{example}
	\begin{example}[Matrix algebras and coordinate-wise derivation]\label{coordinatewise}
		Let $(k,d_k)$ be a differential field.
		The map $\delta^c : M_m(k) \rightarrow M_m(k)$ defined by $\delta^c ((a_{ij})) = (d_k(a_{ij}))$ is a derivation on $M_m(k)$
		and $\delta^c$ is called the \textbf{coordinate-wise derivation on $M_m(k)$}. 
		Then $(M_m(k), \delta^c)$ is a differential algebra over $(k, d_k)$.
	\end{example}
	\begin{example}[Quaternion algebras and standard derivation]\label{DQAstandard}
		Let $(k,d_k)$ be a differential field and $\alpha,\beta\in \mg k$. Let $Q=(\alpha,\beta)_k$ be a quaternion algebra and let $u, v \in Q$ such that $u^2=\alpha$, $v^2=\beta$ and $vu = -uv$. A derivation $d$ on $Q$ such that $d(u)\in k(u)$ and
		$d(v) \in  k(v)$ is called a {\bf standard derivation on $Q$}, and is denoted by $d_{(u,v)}$. 
		This is defined in \cite{AKVRS}.
	\end{example}

	Let $(R,d_R)$ be a differential ring.		
	Elements $c\in R$ such that $d_R(c)=0$ are called {\bf constants.} The kernel of the derivation,  $\{c\in R:d_R(c)=0\}$ forms a ring (indeed a differential ring with the zero derivation).  We call the kernel the {\bf  ring of constants} of $(R,d_R)$ and denote it by $C_{(R,d_R)}$.  A ring homomorphism $\phi: (R,d_R)\to (S,d_S)$ is called a {\bf differential homomorphism} if $d_S(\phi(\alpha))=\phi(d_R(\alpha))$ for all $\alpha\in R$. 
	A ring $S\supseteq R$ is called a  {\bf differential ring extension} of $R$ if there is a differential ring structure on $S$ which is compatible with the differential ring structure on $R$, i.e., $d_S|_R =d_R$. A differential ring extension of a differential field $(k,d_k)$ is a {\bf differential $k$-algebra}. 
	
	The main objects of interest in this article are differential central simple algebras. 
	We have already seen some examples of differential central simple algebras (Examples \ref{coordinatewise}, \ref{DQAstandard}). Classifications of the derivations on certain central simple algebras are noted in \Cref{DCSA}. 
	A crucial aspect of differential central simple algebras is their differential splitting fields.
	These splitting fields can be obtained as  Picard-Vessiot extensions of a matrix differential equation. This connects  differential splitting fields to the differential Galois theory (see \cite{Kol}).      In \Cref{PVE}, we recall the definitions and results related to Picard-Vessiot extensions of a matrix differential equation and their differential Galois groups.   In  \Cref{Dsplitting}, we have  gathered several open questions concerning differential splitting fields. \Cref{tensorprod} explores the tensor products of differential central simple algebras and related properties.
	
	\section{Picard-Vessiot extensions}\label{PVE}
	For  an introduction to differential Galois theory we refer the interested readers to \cite{Kap}, \cite{Kol} and  \cite{Mag94}. 
	For the results in this section,  our main references are \cite[Chapter 1]{PS} and \cite{CHS}.

	Let $(k,d_k)$ be a differential field with its field of constants $C_k$  an algebraically closed field. Let $P\in M_m(k)$. A {\bf Picard-Vessiot ring} $(R,d_R)$ over $(k,d_k)$ for the matrix differential equation $\delta^c(Y)=PY$ is a ring satisfying the following:
	\begin{enumerate}[(a)]
		\item The only ideals $I$ of $R$ satisfying $d_R(I) \subseteq I$  are $(0)$ and $R$.
		\item There exists a matrix $F\in \GL_m(R)$ satisfying $\delta^c(F)=PF$.
		\item The ring $R$ is generated as a ring by $k$, the components of $F$ and the inverse of the determinant of $F$.
	\end{enumerate}
	
	For a matrix differential equation $\delta^c(Y)=PY$ over $(k,d_k)$, there exists a Picard-Vessiot ring and is unique up to $k$-differential isomorphism. 
	A Picard-Vessiot ring is a domain and 
	the fraction field of a Picard-Vessiot ring is called a {\bf Picard-Vessiot extension field} over $k$ and is abbreviated as {\bf PVE} (see \cite[Proposition 1.20, Definition 1.21]{PS}).  
	The next proposition describes the properties of a PVE:
	
	\begin{proposition}\cite[Proposition 1.22]{PS}
		Let $(k,d_k)$ be a differential field and $P\in M_m(k)$.  A differential field extension $(K,d_K)\supset (k,d_k)$ is a PVE for a matrix differential equation $\delta^c(Y)=PY$  if and only if the following hold.
		\begin{enumerate}[(a)]
			\item There exists $F\in \GL_m(K)$ such that $\delta^c(F)=PF$.
			\item The field $K$ is generated over $k$ by the entries of $F$.
			\item The field of constants $C_K=C_k$.
		\end{enumerate}
	\end{proposition}

	Let $(K,d_K)\supset (k,d_k)$ be a PVE for matrix equation $\delta^c(Y)=PY$, where $P\in M_m(k)$. The group of differential $k$-automorphisms of $K$ is called the {\bf differential Galois group} and is denoted by $\dgal(K/k)$. The  differential Galois group can be viewed as a subgroup of $\GL_m(C_k)$ (see \cite[Observations 1.26]{PS}), and hence is a linear algebraic group. The transcendence degree of $K$ over $k$ is equal to the dimension of $\dgal(K/k)$.

	Similar to the classical Galois correspondence, there exists a differential Galois correspondence.
	For a subgroup $H\subset \dgal(K/k)$, let $K^H$ denote the set of elements of $K$ that are fixed by every differential automorphism in $H$. For a differential field $k\subset L\subset K$, let DAut$(L/k)$ denote the set of differential automorphisms in $\dgal(K/k)$ that fixes every element in $L$. We state the differential Galois correspondence in the theorem below.

	\begin{theorem}\cite[Theorem 1]{CHS}\label{Dgalcor}
		For $P\in M_m(k)$, let $(K,d_K)$ be a PVE for a matrix equation $\delta^c(Y)=PY$ over $(k,d_k)$.  Then
		\begin{enumerate}[(a)]
			\item For every (Zariski) closed subgroup $H\subset \dgal(K/k)$, there exists a differential subfield $L\subset K$ containing $k$ such that $K^H=L$ and vice versa.
			\item A differential subfield $L\subset K$ containing $k$ is a PVE over $k$ if and only if $\dgal(K/L)$ is a normal subgroup of $\dgal(K/k)$.
		\end{enumerate}
	\end{theorem}

	\noindent In \cite{CHS}, authors established the Galois
	correspondence for general PVEs, i.e. without the assumption that $C_k$ is algebraically closed.


	\section{Differential central simple algebras}\label{DCSA}
	
	Let $k$ be a field.
	By a {\bf central simple $k$-algebra} we always mean
	a finite-dimensional $k$-algebra that  has no two-sided non-trivial ideal and whose center is $k$. 
	We refer to \cite{Dra},\cite{GilSza} or \cite{BO}, for the standard theory of central simple algebras.
	
	Let $A$ be a central simple $k$-algebra.  
	Let $d_k$ be a derivation on $k$ and $d$ be an extension of $d_k$ to $A$. Then
	$(A,d)$ is a  {\bf differential central simple algebra over $(k, d_k)$}. The set of all derivations on $A$ that are extensions of $d_k$ is denoted by $Der(A/(k,d_k))$.
	Derivations on simple algebras are studied in  \cite{Ami2}, \cite{Hoe} and \cite{Jac}. 
	Initially questions were asked about the existence of extensions of derivations, which were dealt with in \cite{Ho}, \cite{BK},\cite{RS} and  \cite{Ami}.
	
	\begin{theorem}\cite[Theorem, Proposition 1]{Ami} 
		A derivation $d_k$ on the field $k$ can always be extended to a central simple $k$-algebra $A$.
		If $d,d'\in  Der(A/(k,d_k))$, then $d' =d+ \partial_a $ for some $a \in A$.
	\end{theorem}
	
	\noindent In particular, the above result exhibits that one derivation on a central simple $k$-algebra $A$ determines the whole set $Der(A/(k,d_k))$. 
	
	\begin{example}
		Let $d_k$ be the zero derivation on $k$. Then 
		$Der(A/(k,d_k))$ consists of inner derivations on $A$.
	\end{example}
	
	In some cases we can define a special element in $Der(A/(k,d_k))$, which we may call ``standard", depending on the parameter or structure of the algebra.
	In the following example, we see the structure of $Der(M_m(k)/(k, d_k))$ for martix algebras $M_m(k)$.
	
	\begin{example}\label{matrixalgder}
		For the matrix algebra $M_m(k)$, we have $ Der(M_m(k)/(k, d_k)) = \{d_P =\delta^c + \partial_P \mid\ P\in M_m(k), \trace(P) =0\}$ (see \cite[Theorem 2]{Ami}).
	\end{example}

	Symbol algebras are another important example of central simple algebras where a specific derivation is known depending on the defining parameters of the algebra. 
	For $m\geq 2$, assume that $k$ contains a primitive $m$th root of unity  $\zeta$. 
	For  $\alpha, \beta \in \mg k$, the \textbf{symbol algebra} $A = ( \alpha, \beta)_{k, \zeta}$ is an $m^2$-dimensional $k$-algebra generated by $u,v \in A$ satisfying the relations $u^m=\alpha, v^m =\beta$ and $vu = \zeta uv$. 
	Any $4$-dimensional central simple $k$-algebras is isomorphic to a symbol algebra with $m=2$; these are also called  \textbf{quaternion algebras}.  We denote the set of trace $0$ elements in $A$ by $A^0$.
	
	\begin{example}[Symbol algebras and standard derivation] \label{DSA}
		Let $A$ be a $m^2$- dimensional symbol algebra over $k$. A derivation $d\in Der(A/(k,d_k))$ is called a {\bf standard derivation on $A$} if there exists $u,v\in A^0$ such that $u^m, v^m \in \mg k,\  v u=\zeta u v,\ d(u)\in k(u)$ and $d(v)\in k(v)$, and we denote this derivation by $d_{(u,v)}$. For every $d\in Der(A/(k, d_k))$,  there exists a unique element $a \in A^0$ such that $d =d_{(u,v)}+\partial_a$ (see  \cite[Proposition 3.1]{GKS}).
	\end{example}
	
	We now generalize the above definition to an arbitrary central division algebra. Let $D$ be a finite-dimensional central division $k$-algebra. As an easy application of the Double Centraliser Theorem, one gets that $D$ is generated by two elements $D =k[u, v]$ where $u,v$ satisfies certain relations (see \cite[Corollary 3]{HR}). We say $d\in Der(D/(k,d_k))$
	is a {\bf standard derivation on $D$} if there exists $u,v\in D $ such that $D= k[u,v],$ $d(u)\in k(u)$ and $d(v)\in k(v)$, and we denote this derivation by $d_{(u,v)}$.

	In the context of Examples \ref{matrixalgder} and \ref{DSA}, it is natural to ask the following question:
	\begin{question}
		For a central simple algebra $A$ over a differential field $(k,d_k)$, can one completely determine the set of derivations $Der(A/(k, d_k))$ using a standard derivation? 
	\end{question}

	In the context of the above question, one can study the behavior of derivations on subfields of $A$. 
	As we have seen in \Cref{DSA} a standard derivation on a symbol algebra $A$ stabilises two maximal subfields of $A$.
	This motivates the study of subfields of $A$ that are stable under $d$.
	A subfield $L \subseteq A$ is called a \textbf{differential subfield} of $(A,d)$ if $d(L) \subseteq L$.
	Note that, in this case,  $d$ restricts to the unique derivation $d_L$ on $L$. 
	In \cite[Theorem]{Ami}, it was shown that given a central simple $k$-algebra $A$ and a maximal subfield $L$ of $A$, there exists $d\in Der(A/(k,d_k))$ such that $L$ is a differential subfield of $(A,d)$.
	However, there are differential central simple algebras $(A,d)$ having no differential subfields (see \cite{AKVRS} and \cite{GKS}).
	
	\begin{question}
		For a central simple algebra $A$ over  a differential field $(k,d_k)$, classify the set of derivations $d\in Der(A/(k,d_k))$ such that every maximal subfield of $A$ is a differential subfield.
	\end{question}
	
	We say that two differential central simple algebras $(A_1,d_1)$ and $(A_2,d_2)$ over $(k,d_k)$ are {\bf differentially isomorphic} if there exists a  $k$-algebra isomorphism $\phi: A_1\rightarrow A_2$ such that $ \phi \circ d_1 =d_2 \circ\phi$; in this case $\phi$ is called a {\bf differential isomorphism of  $(A_1,d_1)$ and $(A_2,d_2)$}.
	
	\begin{question}
		What is the relation between two  differential central simple algebras having the same differential subfields? Do there exist non-isomorphic differential algebras $(A_1, d_1)$ and $(A_2,d_2)$ over $(k,d_k)$  having same differential subfields?  
	\end{question}

	
	\section{Differential splitting field}\label{Dsplitting}

	Let $A$ be an $m^2$-dimensional central simple $k$-algebra and let $K/k$ be a field extension. 
	Then $A\otimes_k K$ is again a central simple $K$-algebra. We say that {\bf $A$ is split over $K$} if $A\otimes_kK \simeq M_m(K)$.
	It is  well-known  that a central simple algebra splits over its maximal subfields.
	For every subfield $E \subseteq A$, we have $[E:k]$ divides $\sqrt{\dim_k(A)}$.
	Thus, in particular, there exists a field extension $E/k$ with $[E:k] \leq \sqrt{\dim_k(A)}$ and the algebra $A$ is split over $E$.
	We refer the reader to \cite[Chapter 13]{Pierce} for the statements related to maximal subfields and splitting fields.

	Analogously, Juan and Magid \cite{JM} gave a notion of a differential splitting field of a differential central simple algebra. 
	A differential central simple algebra $(A,d)$ over $(k,d_k)$ \textbf{ is split over a differential field $(K,d_K)\supseteq (k,d_k)$} if $(A\otimes_k K, d^{\ast} = d\otimes d_K)$ and  $ (M_m(K),\delta^c)$ are  differentially isomorphic algebras over $(K, d_K)$; we say that $(K, d_K)$ is a \textbf{differential splitting field} id $(A,d)$.
	
	Let $(A,d)$ be a differential central simple algebra over $(k, d_k)$.
	Let $L/k$ be a splitting field of $A$.
	Then there exists a trace zero matrix $P\in M_m(L)$ such that $(A\otimes_kL, d^\ast)\simeq(M_m(L),d_P)$ (see \Cref{matrixalgder}). Therefore,
	in order to understand the splitting of differential central simple algebras, the first step is to understand the splitting of differential matrix algebras. In \cite[Proposition 2]{JM}, a necessary and sufficient condition was obtained for a differential  matrix algebra to be split.
	
	\begin{proposition}\label{JM_split}
		Let $P, Q \in M_m(k)$.  Then differential algebras $(M_m(k), d_P)$ and $(M_m(k), d_Q)$ are differentially isomorphic if and only if there exists $F \in GL_m(k)$ such that $F^{-1}\delta^c(F)+F^{-1}QF =P$.
		In particular, $(M_m(k), d_P)$ is split if and only if there exists $F \in GL_m(k)$ such that $\delta^c(F) = FP$. 
	\end{proposition}
	\noindent Note that, for  $F^{-1}$ we have $\delta^c(F^{-1}) = PF^{-1}$ and hence we can as well say that $(M_m(k), d_P)$ is split if and only if there exists $F \in GL_m(k)$ such that $\delta^c(F) = PF$. 
	This leads to the following result for arbitrary differential central simple algebras.
	\begin{theorem}\cite[Corollary 1]{JM}
		A differential central simple algebra  over $(k, d_k)$ splits  over a PVE of $(k,d_k)$.
	\end{theorem}
	
	Differential splitting fields of a differential central simple algebra are not necessarily algebraic extensions (see \cite[Example 4.3, Theorem 4.5]{AKVRS}, \cite[Theorem 5.5]{GKS} and \cite[Example 5.2, Theorem 6.5]{GKS2}). It is intriguing to determine the transcendence degree of a differential splitting field for a given differential central simple algebra. We set
	$$\ddeg_{k}(A,d) := \min\{ \trdeg_k(E)\mid (E,d_E)\supseteq (k,d_k) \mbox{ splits } (A,d) \}$$
	and call $\ddeg_{k}(A,d)$ the {\bf differential splitting degree of $(A,d)$}.
	
	It is clear form \Cref{JM_split} that $\ddeg_{k}(A,d)$ is bounded by $\dim_k(A)$. Let $(K,d_K)$ be a PVE of $(k,d_k)$ splitting $(A,d)$ and let $\dgal(K/k)$ be its differential Galois group.
	Since we can choose the matrix $P$ to be a trace zero matrix, $\dgal(K/k)$ can be viewed as a subgroup of $\SL_m(C_k)$, by \cite[Exercise 1.35 (5a)]{PS}. 
	Then   
	$$\trdeg_k(K) = \dim_k(\dgal(K/k))\leq m^2-1.$$
	Thus $\ddeg_{k}(A,d)$ is bounded by $m^2-1$.

	\begin{question}
		Let $A$ be a central simple $k$-algebra. 
		Which numbers $i$ between $0\leq i\leq m^2-1$ can occur as $\ddeg_{k}(A,d)$ for derivations $d\in Der(A/(k,d_k))$?  \end{question}
	
	Examples of differential quaternion algebras $(Q,d)$ over $(k, d_k)$ such that  $\ddeg_{k}(Q,d)=i$ for each $0\leq i\leq 3$ were produced in \cite{GKS2}.
	
	\begin{example}\cite[Example 5.2 and Theorem 6.5]{GKS2}
		Let $(k=\mathbb Q(t),d_k)$ be a differential field  with $d_k(t)=1$. Let $Q=(1,t)_k$ be a quaternion algebra and $u,v\in Q^0$ such that $u^2=1$, $v^2=t$ and $vu=-uv$. Consider a derivation $d=d_{(u,v)}+\partial_a$ on $Q$, where $a\in Q$. Then
		\begin{enumerate}[(a)]
			\item  For $a =0$, $\ddeg_k(Q,d)=0$.
			\item  For $a =\frac{1}{t}v$, $\ddeg_k(Q,d)=1$.
			\item   For $a =\frac{1}{4t}(-u-2v+2uv)$, $\ddeg_k(Q,d)=2$.
			\item  For $a =-\frac{1}{4t}u-v$, $\ddeg_k(Q,d)=3$.
		\end{enumerate}
	\end{example}
	Using \cite[Theorem 6.5]{GKS2}, one can produce several examples of differential quaternion algebras with  $\ddeg_k(Q,d)=0,1$. On the other hand, for $\ddeg_k(Q,d)=2,3$ we only know examples of split quaternion algebras. Therefore, we ask the following question.
	\begin{question}
		Does there exist a differential division quaternion algebra $(Q,d)$ such that $\ddeg_k(Q,d)$ equals $2$ or $3$?
	\end{question}
	
	Since differential Galois groups of a PVE over $k$  can be viewed as a subgroup of $\SL_m(C_k)$, there may exist a relation between differential splitting fields of differential central simple algebras and subgroups of $\SL_m(C_k)$. 
	\begin{question}
		Which linear algebraic subgroups of  $\SL_m( C_k)$ can occur as Galois groups of PVE over $k$ that splits a differential central simple algebra over $k$?  
	\end{question}

	It is also interesting to relate differential subfields to differential splitting fields
	For example, it was shown in \cite[Theorem 5.4]{GKS}, that a symbol algebra $(\alpha, \beta)_{k, \zeta}$ with a standard derivation (where we have differential subfields that are also maximal) is split by a finite extension of $k$.

	\begin{question}
		Let $(L,d_L)$ be a PVE of $(k,d_k)$ splitting $(A,d)$.  Explore the relation between differential subfields of $(A,d)$ and $(L,d_L)$. More precisely we can ask:
		Assuming that there is a maximal subfield $M \subseteq A$ which is also a differential subfield of $A$, what can we say about $\ddeg_k(A,d)$? 
		
	\end{question}

	\begin{question}
		Let $A$ be a central simple $k$-algebra and let $d_1, d_2 \in \Der(A/(k,d_k))$. 
		Asssume that $(A, d_1)$ and $(A, d_2)$ have the same differential splitting fields. Is $(A, d_1) \simeq (A, d_2)$?
	\end{question}
	
	The following characterises differential quaternion algebras over $(k, d_k)$ that are split over a finite differential field extension of $(k, d_k)$.
	
	\begin{theorem}\cite[Theorem 7.3]{GKS2}\label{finitesplit}
		A differential quaternion algebra  splits over  a finite field extension of $k$ if and only if the derivation on  the quaternion algebra becomes standard over a finite extension of $k$.
	\end{theorem}

	\begin{question}
		Does \Cref{finitesplit} hold for symbol algebras and division algebras? 
	\end{question}

	\section{Tensor products of differential central simple algebras}\label{tensorprod}
	
	Given two differential central simple algebras $(A_1,d_1)$ and $(A_2,d_2)$ over $(k,d_k)$, we obtain the tensor differential algebra as follows.
	The algebra $A_1\otimes_k A_2$ is again a  central simple $k$-algebra and the map $d_1\otimes d_2 : A_1\otimes_k A_2 \rightarrow A_1\otimes_k A_2 $ defined by  $(d_1\otimes d_2)(a_1\otimes a_2) =  d_1(a_1)\otimes a_2 + a_1\otimes d_2(a_2)$ is a derivation on $A_1\otimes_k A_2$. 
	We write $(A_1,d_1) \otimes_{(k,d_k)} (A_1,d_1)$ for the  differential central simple algebra $(A_1\otimes_k A_2, d_1\otimes d_2)$ over $(k,d_k)$. 
	
	\begin{question}
		Let $A$ be a central simple $k$-algebra and let $d_1, d_2 \in \Der(A/(k,d_k))$. 
		Asssume that $(A, d_1)\otimes_k (A, d_2)$ is split over $(k,d_k)$.
		Are $(A,d_1)$ and $(A,d_2)$ differentially isomorphic?
	\end{question}
	
	The Brauer group of a field $k$, denoted by $\br k$, is   the quotient group of  the group generated by  central simple algebras modulo the subgroup generated by matrix algebras. Due to Weddurburn's Theorem, every  central simple $k$-algebra  $A$ is Brauer equivalent to a unique division algebra up to isomorphism.

	Analogous to the usual definition of the Brauer group of a field, a definition of the differential Brauer group was given in \cite{Hoo}. 
	However, with the definition in \cite{Hoo} it turned out that when the base field is of characteristic zero the differential Brauer group is the same as the usual Brauer group and thus does not contain any extra information. 
	In \cite{JM}, the {\bf differential Brauer group} denoted by $Br^{diff}(k,d_k)$ is defined as the quotient group of  the group generated by differential matrix algebras over $(k,d_k)$ modulo the subgroup generated by matrix algebras with coordinate-wise derivation. Note that $Br^{diff}(k,d_k)$ encompasses information only about differential matrix algebras. 
	When $k$ is algebraically closed (or a $C_1$-field) then all central simple algebras are matrix algebras.

	It is  well-known that $\br k$ is a torsion group, 
	that is every element in $\br{k}$ has finite order.  
	It is well-known for a central division $k$-algebra $D$ that the order of the class of $D$ in $\br{k}$ divides $\sqrt{\dim_k(D)}$.  
	However, it is possible to construct a differential matrix algebra $(M_m(k),d)$ such that for every $n\in \nat $, $(M_m(k),d) \otimes_k \cdots \otimes_k (M_m(k),d)$ ($n$-times) in a non-split differential algebra.
	A differential central simple algebra $(A,d)$ over $(k,d_k)$ is said to be of \textbf{finite order} if there exists an  integer $n\in \nat $ such that $(A,d) \otimes_k \cdots \otimes_k (A,d)$ ($n$-times) is split differential  algebra over $(k,d_k)$.
	The following question is motivated by the work of K. Singla and A. Kulshrestha \cite[Theorem 5.1]{AKKS}, where for certain  matrix differential algebras,  they have shown that finite order implies splitting by a finite extension, i.e.~ $\ddeg =0$.   
	
	\begin{question}
		Let $(A,d)$ be a differential central simple algebra over $(k,d_k)$ of finite order. Does it imply that $\ddeg(A,d) =0$? 
	\end{question}
	
	Since every element in $\br{k}$ has finite order, it is enough to consider the above question for matrix algebras. One can also try to deal with the case of quaternion algebras or symbol algebras separately.

	\begin{question}
		Find relation between differential splitting degrees of two Brauer equivalent differential central simple algebras.
	\end{question}

	\subsection*{Acknowledgement}
	Authors wish to express their gratitude to Varadharaj R.~Srinivasan and  Amit Kulshrestha for  several insightful discussions on this topic which has given rise to many questions in this article.
	We also thank  Kanika Singla for discussions happened on this topic during a seminar series organised by Varadharaj R.~Srinivasan. 
	We also thank the referee for valuable comments to improve the manuscript.

	\bibliographystyle{amsalpha}

\end{document}